\documentclass[12pt]{amsart}
\usepackage{amsfonts, amsmath, amsthm}
\usepackage{epsfig}
\usepackage{psfrag}

\newtheorem{pro}{Proposition}[section]
\newtheorem{thm}[pro]{Theorem}
\newtheorem{lem}[pro]{Lemma}
\newtheorem{clm}[pro]{Claim}

\newtheorem{cor}[pro]{Corollary}

\newtheorem*{claim}{Claim}

\theoremstyle{definition}
\newtheorem{dfn}[pro]{Definition}

\theoremstyle{remark}

\title{Distance and bridge position} 
\date{\today}
\address{D. Bachman, Mathematics Department, California Polytechnic State University}
\email{dbachman@calpoly.edu}
\address{S. Schleimer, Mathematics Department, University of Illinois at Chicago}
\email{saul@math.uic.edu}
\author{David Bachman}
\author{Saul Schleimer}
\begin{document}
\begin{abstract}
J. Hempel's definition of the {\it distance} of a Heegaard surface generalizes to a complexity for a knot which is in bridge position with respect to a Heegaard surface. Our main result is that the distance of a knot in bridge position is bounded above by twice the genus, plus the number of boundary components, of an essential surface in the knot complement. As a consequence knots constructed via sufficiently high powers of pseudo-Anosov maps have minimal bridge presentations which are thin. 
\end{abstract}
\maketitle

\noindent
Keywords: Heegaard Splitting, Curve Complex

\section{Introduction.}
Hempel's definition of the {\it distance} of a Heegaard splitting \cite{hempel:01} is a natural measure of complexity, generalizing the standard notions of {\it reducibility} (distance zero), {\it weak reducibility} (distance at most one), and {\it strong irreducibility} (distance at least two). Hempel proves that there exist Heegaard splittings of arbitrarily high distance. 

In his Ph.D. thesis, K. Hartshorn related the distance of a Heegaard splitting to the genus of any essential surface:

\medskip
\noindent {\bf Theorem} \cite{hartshorn:99}. {\it Let $M$ be a closed, orientable, irreducible 3-manifold with Heegaard splitting $F$. Suppose $M$ contains an orientable, incompressible surface $S$. Then the distance of $F$ is bounded above by twice the genus of $S$.}
\medskip

\noindent Hartshorn's result refines work of T. Kobayashi \cite{Kobayashi88}.

We begin this paper by recalling a generalization of the curve complex for surfaces with non-empty boundary. This allows us to translate Hempel's definition of distance for Heegaard splittings to a definition of distance for knots which are in bridge position with respect to a Heegaard surface \cite{ms:91}. Our main result is a translation of Hartshorn's Theorem into this new context:

\medskip
\noindent {\bf Theorem \ref{t:main}.} {\it Let $K$ be a knot in a closed, orientable 3-manifold $M$ which is in bridge position with respect to a Heegaard surface $F$. Let $S$ be a properly embedded, orientable, essential surface in $M_K$. Then the distance of $K$ with respect to $F$ is bounded above by twice the genus of $S$ plus $|\partial S|$.}
\medskip

In the special case of a meridional disk we find that a stronger result holds; the distance of $K$ with respect to $F$ is zero. This follows from a variant of the Haken Lemma \cite{haken:68} (see Lemma \ref{l:Haken}).

Although our proof contains Hartshorn's result as a special case ($K=\emptyset$), there is an interesting qualitative difference. Unlike Hartshorn, we make no minimality assumption on the way in which $S$ intersects $F$. That is, {\it any generic position} of $S$ with respect to $F$ forces the bound on distance as stated in the theorem. 

The main idea behind our proof is to simply count saddles. Let $d(K,F)$ denote the distance of $K$ with respect to $F$. It is a standard technique in 3-manifold topology to use a Heegaard splitting $F$ for a 3-manifold $M$ to define a height function $h$ on $M$. This, in turn, induces a height function on a surface $S$ in $M$. With respect to this height function $S$ will have maxima, minima, and saddles. The ``moral" is that each critical point of $S$ either 

\begin{enumerate}
    \item contributes at most 1 to $d(K,F)$ and exactly -1 to the Euler characteristic of $S$ or
    \item contributes nothing to $d(K,F)$ and nothing to the Euler characteristic of $S$.
\end{enumerate}

Hence, the distance of $K$ with respect to $F$ would then be bounded by the negative of the Euler characteristic of $S$. Unfortunately, for Heegaard splittings the above classification isn't exactly correct. We find that there may be at most two ``special" critical points that each contribute one to the distance of $K$, but nothing to the Euler characteristic of $S$. This gives the bound $d(K,F) \le -\chi(S) +2 = 2g(S)+|\partial S|$. 

We note that several authors have explicitly computed the distances of various classes of knots (using varying definitions of {\it distance}). See, for example, \cite{sakuma:99}, \cite{morimoto:89} and \cite{saito:03}. 

In the final section we present corollaries to Theorem \ref{t:main}. Among these are:

\medskip
\noindent {\bf Corollary \ref{c:c1}.} {\it Suppose $K$ is a knot in $S^3$ whose distance is $d(K,F)$ with respect to a bridge sphere $F$. Then the genus of $K$ is at least $\frac {d(K,F)-1}{2}$.}

\medskip

\noindent {\bf Corollary \ref{c:c2}.} {\it If $K$ is a knot whose distance is at least 3 with respect to some Heegaard surface, then the complement of $K$ is hyperbolic.}

\medskip

Finally, we define the {\it bridge link} associated to an element of the braid group $B_{2n}$ to be the link obtained by gluing two trivial $n$-strand tangles by this element. By a construction essentially due to Kobayashi \cite{Kobayashi88} powers of certain pseudo-Anosov maps give associated bridge links with arbitrarily high distance. Suppose $\phi$ is such a map. Then it follows from Corollary \ref{c:thinbridge} that for all sufficiently high powers of $\phi$ if the associated link is a knot then its minimal bridge presentation is thin.

It is possible, {\it a priori}, that bridge knots associated to high powers of pseudo-Anosov maps have low bridge numbers. The following conjecture would rule this out:

\medskip

\noindent {\bf Conjecture.} {\it Suppose $K$ is a knot whose distance is at least 2 with respect to some Heegaard surface $F$. Then the distance of $K$ with respect to any other Heegaard surface is bounded above by $\chi(F-K)+2$.}

\medskip

Compare this to the statement of Theorem \ref{t:main}. In the theorem we assert that the distance of a knot with respect to a Heegaard surface is bounded by two plus the Euler characteristic of an essential surface. In the conjecture we propose that distance is similarly bounded by a strongly irreducible surface.

\section{Basic Definitions.}
In this section we give the definitions that will be used throughout the paper. Let $K$ be a knot in a closed, orientable  3-manifold, $M$. Let $M_K=M-N(K)$ where $N(K)$ denotes a regular neighborhood of $K$. For the remainder of this paper all surfaces $S$ in $M_K$ will be embedded, compact, and orientable with $S \cap \partial M_K \subset \partial S$.

\begin{figure}[htbp]
        \psfrag{g}{$\gamma$}
        \vspace{0 in}
        \begin{center}
        \epsfxsize=4 in
        \epsfbox{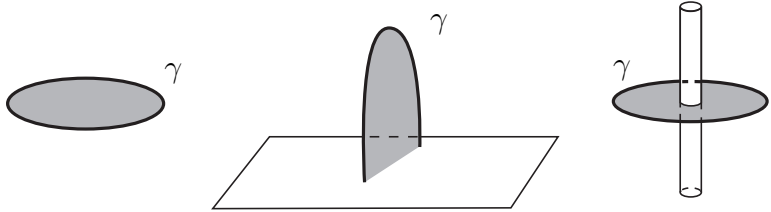}
        \caption{Disk, bigon, and meridional cut surfaces.}
        \label{f:cutsurface}
        \end{center}
        \end{figure}

\begin{dfn}
A {\it cut surface} (see Figure \ref{f:cutsurface}) is either
\begin{enumerate}
    \item a disk $E \subset M_K$ such that $E \cap \partial M_K=\emptyset$,  
    \item a bigon $E \subset M_K$ such that $E \cap \partial M_K$ is an arc, or
    \item an annulus $E \subset M_K$ with exactly one meridional boundary component on $\partial M_K$. In other words, $E \cap \partial M_K$ is a loop bounding a disk in $\overline{N(K)}$.
\end{enumerate}
If $E$ is a cut surface and $\gamma =\overline{\partial E-\partial M_K}$ then we say $\gamma$ {\it bounds a cut surface}. 
\end{dfn}

A properly embedded simple curve in $S$ is {\it inessential} if it bounds a subsurface of $S$ which is a cut surface, and {\it essential} otherwise.

Suppose $\gamma$ bounds a cut surface $E$, $S$ is properly embedded in $M_K$, and $S \cap E=\gamma$. We may then {\it surger} $S$ along $E$ by replacing a neighborhood of $\gamma$ in $S$ with two parallel copies of $E$. If $\gamma$ is essential in $S$ then we say $E$ is a {\it compression} for $S$. In this case we also say $\gamma$ {\it bounds a compression} for $S$. 

A properly embedded surface $S \subset M_K$ is {\it essential} if first there are no curves on $S$ which bound compressions in $M_K$ and second $\partial S$ (if non-empty) is not null-homotopic on $\partial M_K$. We also consider a  2-sphere to be essential if it does not bound a ball in $M_K$. Note that this notion of essential is {\it not} identical to
that of ``super-incompressible" found in \cite{MorganBass84}.

A {\it handlebody} is a 3-manifold which is homeomorphic to the closure of a regular neighborhood of a compact, connected graph in ${\mathbb R}^3$. If such a graph has no valence one vertices, and the corresponding handlebody has non-zero genus, then the graph's image under such a homeomorphism is a {\it spine} of the handlebody. We will insist that the spine of a 3-ball be a single edge.

A closed surface $F$ in $M$ is a {\it Heegaard surface} of M if $F$ separates $M$ into two handlebodies. An arc  properly embedded in $H$ is {\it trivial} if it bounds a bigon in $H$. Suppose $K$ is a knot in a 3-manifold $M$ with Heegaard surface $F$. The knot $K$ is in {\it bridge position with respect to $F$} \cite{ms:91} if $K$ meets each of the handlebodies bounded by $F$ in a collection of trivial arcs. Such a position is sometimes referred to as a {\it $(g,b)$-presentation} of $K$, where $g=\mbox{genus}(F)$ and $2b=|K \cap F|$.

\section{The Arc Complex}

Following Hempel's definition of the distance of a Heegaard splitting \cite{hempel:01} we now define the {\it distance} of a knot $K$ which is in bridge position with respect to a Heegaard surface $F \subset M$. Let $M_K=M-N(K)$ and let $F_K=F \cap M_K$.

Construct a 1-complex $\Gamma (F_K)$ as follows: for each proper isotopy class of essential curve in $F_K$ there is a vertex of $\Gamma (F_K)$. There is an edge of $\Gamma (F_K)$ between two distinct vertices if and only if there are representatives of the corresponding isotopy classes which are disjoint. $\Gamma (F_K)$ is called the {\it arc complex} of $F_K$ (see, for example, \cite{mm:99}).

Now, $F_K$ divides $M_K$ into two submanifolds, $H$ and $H'$. Let $V$ and $V'$ denote the sets of vertices of $\Gamma (F_K)$ which correspond to curves which bound compressions in $H$ and $H'$ (resp.). Then $d(K,F)$, the {\it distance of $K$ with respect to $F$}, is defined to be the number of edges in the shortest path from $V$ to $V'$ in $\Gamma (F_K)$. As long as $\chi(F_K) \leq -2$ this is well defined, since the arc complex is connected in those cases.

%\begin{lem}
%\label{l:MeridionalCompression}
%Let $K$ be a knot in a 3-manifold $M$. If there is a meridional disk in $\partial M_K$ then either $M$ contains an essential 2-sphere or $K=\{{\rm pt}\} \times S^1 \subset M= S^2 \times S^1$. 
%\end{lem}

\section{Lemmas}

The following is a slight variant of the Haken Lemma \cite{haken:68}. We assume that the reader is familiar with W. Jaco's proof of this result (Theorem II.7 of \cite{jaco:80}). 

\begin{lem}[Haken]
\label{l:Haken}
Let $K$ be a knot in a 3-manifold $M$ which is in bridge position with respect to a Heegaard surface $F$. If $M_K$ contains an essential 2-sphere or meridional disk then $d(K,F)=0$.
\end{lem}

\begin{proof}
Among all essential 2-spheres and  meridional disks in $M_K$ choose one, $S$, meeting $F_K$ minimally. Let $H$ and $H'$ denote the submanifolds of $M_K$ bounded by $F_K$, with $\partial S$ (if nonempty) contained in $H$. If $S \cap F_K=\emptyset$ then $S$ lies entirely in $H$ or $H'$, a contradiction. It follows that $S \cap F_K$ is a non-empty set of loops which are essential on $F_K$. Hence, if $S$ meets $F_K$ in a single loop then the result follows. 

Suppose then that $|S \cap F_K|>1$. Let $H^*$ denote one of $H$ or $H'$, where there is a component $T$ of $S \cap H^*$ with $|\partial T-\partial S|\ge 2$. Choose a {\it basis} $\Lambda$ for $H^*$, {\it i.e.} a system of disks and bigons which cuts $H^*$ into a 3-ball. If $S \cap \Lambda$ contains any loops then surger $S$ along these, innermost (on $\Lambda$) first. Note that at least one component of the resulting surface is again an essential sphere or meridional disk. We continue to denote this surface by $S$. 

We now reduce $|S \cap \Lambda|$ via the following procedure. If any component of $(S \cap H^*)-\Lambda$ is a bigon then surger $\Lambda$ along this surface. Some subcollection of the resulting set is again a basis, which we continue to denote by $\Lambda$. If not then choose a bigon of $\Lambda -S$, and use this to guide an isotopy of $S$ (see the ``isotopy of type A" on page 24 of \cite{jaco:80}). Repeat this procedure until all components $T$ of $S \cap H^*$ have $|\partial T-\partial S|=1$. Let $S'$ denote the resulting surface.  

It follows from the argument of Lemma II.9 of \cite{jaco:80} that if $H^*=H'$ then $|S' \cap F_K|<|S \cap F_K|$ and we have reached a contradiction. If $H^*=H$ then $|S' \cap F_K| \le |S \cap F_K|$. If equality holds we repeat the above with $S'$ replacing $S$ and letting $H^*=H'$. This gives a surface $S''$ with $|S'' \cap F_K|<|S \cap F_K|$, a contradiction. 
\end{proof}

\begin{lem}
\label{l:NoTypeMismatch}
Let $K$ be a knot in a 3-manifold $M$ which is in bridge position with respect to a Heegaard surface $F$. If $\gamma$ bounds two cut surfaces $A$ and $B$ with $A \cap B=\gamma$ then either both are bigons, both annuli, both disks, or $d(K,F)=0$.
\end{lem}

\begin{proof}
If $A$ and $B$ are different types then their union is a meridional disk. The result now follows from Lemma \ref{l:Haken}.
\end{proof}

\begin{lem}
\label{l:CutImpliesCompression}
Let $K$ be a knot in a 3-manifold $M$ which is in bridge position with respect to a Heegaard surface $F$ and let $Q$ be any properly embedded surface in $M_K$. If there is a curve $\gamma$ which is essential on $Q$ and bounds a cut surface $E$ in $M_K$ then either there is a curve $\gamma' \subset E \cap Q$ which bounds a compression for $Q$, or $d(K,F)=0$.
\end{lem}

\begin{proof}
Let $\Lambda \subset E \cap Q$ be the collection of curves which are essential on $Q$. Let $E'$ denote the closure of a component of $E-\Lambda$ which is a cut surface. Let $\gamma'=E' \cap \Lambda$. Consider the set $\Theta$ of cut surfaces bounded by $\gamma'$ such that the only curve of intersection with $Q$, essential on $Q$, is $\gamma'$. Note that $E'$ is such a surface, so $\Theta$ is non-empty. Let $E^*$ be an element of $\Theta$ with $|E^* \cap Q|$ minimal. 

We now claim  $E^* \cap Q=\gamma'$. Suppose not. Let $E''$ be a cut surface component of $E^*-Q$. The curve $\gamma''=E'' \cap Q$ is inessential on $Q$ and hence bounds two cut surfaces, $A \subset Q$ and $E''$. Note that $A \cap E'' = \gamma''$. By Lemma \ref{l:NoTypeMismatch} we may obtain a new cut surface from $E^*$ by replacing $E''$ with a push-off of $A$. This violates the minimality of $|E^* \cap Q|$. We conclude $E^*$ is a compression for $Q$, which finishes the proof. 
\end{proof}

\begin{lem}
\label{l:EssentialPersistsCutDisk}
Let $K$ be a knot in a 3-manifold $M$ which is in bridge position with respect to a Heegaard surface $F$ and $S$ an essential surface in $M_K$. If we surger $S$ along a disk or bigon cut surface then at least one of the remaining components is essential, or $d(K,F)=0$.
\end{lem}

\begin{proof}
By assumption there is a curve $\gamma \subset S$ which bounds a cut surface $E'$, homeomorphic to a disk, such that $E' \cap S=\gamma$. Since $S$ is essential, $\gamma$ bounds a cut surface $E \subset S$. Surgering $S$ along $E'$ then produces two surfaces, isotopic to $E \cup E'$ and $S'=(S-E) \cup E'$. Suppose $S'$ is not essential. Let $\gamma'$ bound a compression $C$ for $S'$. As $E'$ is homeomorphic to a disk we may properly isotope $\gamma'$ off of $E'$. The curve $\gamma'$ is now on $S$, and bounds the cut surface $C$. By Lemma \ref{l:CutImpliesCompression} there is a compression $C'$ for $S$, a contradiction.
\end{proof}

\begin{lem}
\label{l:EssentialPersists}
Let $K$ be a knot in a 3-manifold $M$ which is in bridge position with respect to a Heegaard surface $F$ and $S$ an essential surface in $M_K$. If we surger $S$ along a cut surface then at least one of the remaining components is essential, or $d(K,F)=0$.
\end{lem}

\begin{proof}
By assumption there is a curve $\gamma \subset S$ which bounds a cut surface $E'$ such that $E' \cap S=\gamma$. Since $S$ is essential, $\gamma$ bounds a cut surface $E$ in $S$. Surgering $S$ along $E'$ then produces two surfaces, isotopic to $E \cup E'$ and $S'=(S-E) \cup E'$. 

By Lemma \ref{l:EssentialPersistsCutDisk} we may assume $E'$ is an annulus. By Lemma \ref{l:NoTypeMismatch} we may assume $E$ is also an annulus. If $E \cup E'$ is essential then we are done. Otherwise there must be a compressing bigon $B$ for $E \cup E'$ (since the core loop of $E \cup E'$ is not essential). Surgering $E \cup E'$ along $B$ gives a disk $D$ with $\partial D \subset \partial M_K$ bounding a disk $D' \subset \partial M_K$. If the sphere $D \cup D'$ is essential then the proof is complete by Lemma \ref{l:Haken}. Otherwise we conclude $E \cup E'$, together with an annulus of $\partial M_K$, bounds a solid torus. If the interior of the solid torus is disjoint from $S$ then $S'$ is properly isotopic to $S$ and we are done. If $S$ meets the interior of the solid torus then we may push it entirely into the solid torus. Now consider $B \cap S$. Some component of $B-S$ is then a cut surface for $S$. This cut surface is either a disk or a bigon. By Lemma \ref{l:EssentialPersistsCutDisk} we may surger $S$ along this cut surface and obtain another essential surface which meets $B$ fewer times. Continuing in this way we obtain an essential surface inside the solid torus which misses $B$, and is hence contained in a ball. This is impossible.
\end{proof}

\section{Proof of Main Theorem}
\label{s:main}

We recall the statement. 

\begin{thm}
\label{t:main}
Let $K$ be a knot in a closed, orientable 3-manifold $M$ which is in bridge position with respect to a Heegaard surface $F$. Let $S$ be a properly embedded, orientable, essential surface in $M_K$.  Then the distance of $K$ with respect to $F$ is bounded above by twice the genus of $S$ plus $|\partial S|$.
\end{thm}

We now begin the proof. Throughout we assume that $d(K,F)>0$ to avoid the special cases of Lemmas \ref{l:Haken} through \ref{l:EssentialPersists}. Let $\Sigma _0$ and $\Sigma _1$ denote spines of the handlebodies bounded by $F$. Let $h:M \rightarrow I$ denote a height function on $M$, so that $h^{-1}(0)=\Sigma _0$ and $h^{-1}(1)=\Sigma _1$. Also we require that for every $t \in (0,1)$ the surface $h^{-1}(t)$ is parallel to $F=h^{-1}(1/2)$. Because $K$ is in bridge position with respect to $F$ we may isotope $K$ so that each arc of $K -F$ has one critical point with respect to $h$. Now pull each minimum down to $\Sigma _0$ and each maximum up to $\Sigma_1$. Note that if $M=S^3$ and $F$ is a sphere then we may assume that $K$ has at least two maxima and at least two minima. In this case $\Sigma_0$ and $\Sigma _1$ are edges, and we assume that the vertices $\partial \Sigma _0$ coincide with two minima of $K$ and the vertices $\partial \Sigma _1$ coincide with two maxima.

Let $F(t)=h^{-1}(t) \cap M_K$. Let $H(t)$ be the closure of the component of $M_K-F(t)$ which contains $\Sigma _0$. Let $H'(t)$ be the closure of $M_K-H(t)$. Let $\epsilon _0$ be chosen just larger than the radius of $N(K)$, but small enough so that $S$ meets $H(\epsilon _0)$ and $H'(1-\epsilon _0)$ in compressions for $F(\epsilon_0)$ and $F(1-\epsilon_0)$. Then the surface $F(t)$ is homeomorphic to $F_K=F \cap M_K$ for every value of $t \in [\epsilon_0,1-\epsilon_0]$. Hence, the submanifold $\bigcup _{t=\epsilon_0}^{1-\epsilon_0} F(t)$ is homeomorphic to $F_K \times [\epsilon_0,1-\epsilon_0]$. Let $\pi$ denote the composition of such a homeomorphism with projection onto the first factor. Hence, if $\gamma$ is a curve on $F(t)$, for some $t\in [\epsilon_0,1-\epsilon_0]$, then $\pi(\gamma)$ is a curve on $F_K$.

We make two types of assumptions on the position of the essential surface $S$. Any surface whose position satisfies these assumptions we will say is in {\it standard position}.  The first concerns how $S$ meets $\partial M_K$ and the second is a genericity assumption on the interior of $S$. Near the boundary of $S$ we assume the following:
\begin{itemize}
    \item Meridional boundary components are ``level." That is, if $S$ has meridional boundary then for each boundary component $C$ of $S$ there exists a $t \in (\epsilon_0,1-\epsilon_0)$ such that $C \subset \partial F(t)$. We consider $t$ a critical value for $S$ if some boundary component of $S$ is contained in $\partial F(t)$. 
    \item If $S$ does not have meridional boundary then for generic $t$ and each component $\gamma$ of $\partial S-F(t)$ the endpoints of $\gamma$ lie on distinct boundary components of $F(t)$. 
\end{itemize}

\begin{figure}[htbp]
        \psfrag{F}{$F(t+\epsilon)$}
        \psfrag{G}{$F(t-\epsilon)$}
        \psfrag{M}{$\partial M_K$}
        \psfrag{S}{$P$}
        \psfrag{C}{$C$}
        \vspace{0 in}
        \begin{center}
        \epsfxsize=4.5 in
        \epsfbox{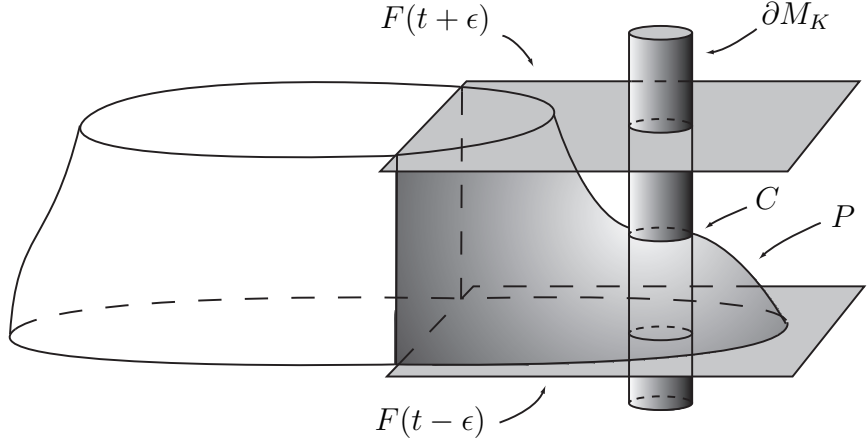}
        \caption{A piece of $S$ between levels, $F(t-\epsilon)$ and $F(t+\epsilon)$, before and after a meridional boundary component, $C$.}
        \label{f:annulus}
        \end{center}
        \end{figure}

Note that the above is possible since $\partial S$ is not null-homotopic on $\partial M_K$. In the interior of $M_K$ we assume the position of $S$ is generic in the following sense: 
\begin{itemize}
    \item All critical points of $h|_S$ are maxima, minima, or saddles. We will refer to any such critical point whose height is between $\epsilon_0$ and $1-\epsilon _0$, and any meridional boundary component, as a {\it critical submanifold} (of $S$). 
    \item The heights of any two critical submanifolds of $S$ are distinct.
    \item Suppose a meridional boundary component $C$ of $S$ happens at height $t$. Let $P$ denote the closure of the component of $S-F(t \pm \epsilon)$ that has $C$ as a boundary component. Then $P$ is a once-punctured annulus with one boundary component on each of $F(t-\epsilon)$ and $F(t+\epsilon)$ (see Figure \ref{f:annulus}). (This uses the fact that $\partial M_K$ is connected.)
\end{itemize} 

\begin{clm}
\label{c:noessentialinhandlebody}
For each $t \in [\epsilon_0,1-\epsilon_0]$ the submanifolds $H(t)$ and $H'(t)$ of $M_K$ do not contain any essential surfaces.
\end{clm}

\begin{proof}
Choose a basis $\Lambda$ of compressing disks and bigons in $H(t)$ which cut it into a ball. Suppose $D \in \Lambda$. Let $D'$ be a cut surface component of $D-Q$, where $Q$ is some essential surface in $H(t)$. By Lemma \ref{l:EssentialPersistsCutDisk} compressing $Q$ along $D'$ yields an essential surface which meets $D$ fewer times. Continuing in this way we produce an essential surface in $H(t)$ disjoint from $\Lambda$, and hence in a ball. This is impossible.
\end{proof}

\begin{dfn}
Let $t_0=\sup \{t \in [\epsilon_0,1-\epsilon_0]\thinspace \big{|}$ some curve in $S \cap F(t)$ bounds a compression for $F(t)$ in $H(t)\}$. Note that the compression for $F(t)$ need not be a subsurface of $S$. Let $t_1=\inf \{t \in [t_0,1-\epsilon_0]\thinspace \big{|}$ some curve in $S \cap F(t)$ bounds a compression for $F(t)$ in $H'(t)\}$. 
\end{dfn}

\begin{clm}
The values $t_0$ and $t_1$ are well defined, and $t_0 > \epsilon_0$. 
\end{clm}

\begin{proof}
To establish the claim we show that for some small $\epsilon> \epsilon_0$ there are curves in  $S \cap F(\epsilon)$ and $S \cap F(1-\epsilon)$ which bound compressions for $F(\epsilon)$ and $F(1-\epsilon)$ (resp.) in $H(\epsilon)$ and $H'(1-\epsilon)$ (resp.). Hence, $t_0$ and $t_1$ are well defined. 

There are essentially two cases. Suppose first the essential surface $S$ is closed, or has meridional boundary. If $S \cap \Sigma _0=\emptyset$ then $S$ can be properly isotoped entirely into $H'(\epsilon)$, violating Claim \ref{c:noessentialinhandlebody}. We conclude that $S \cap \Sigma _0 \ne \emptyset$. $F(\epsilon) \cap S$ then contains a loop which bounds a compression for $F(\epsilon)$ in $H(\epsilon)$. On the other hand, if $S$ has non-empty, non-meridional boundary then $F(\epsilon) \cap S$ contains an arc which bounds a bigon compression in $H(\epsilon)$. This proves $t_0$ is well defined and $t_0>\epsilon>\epsilon _0$. A symmetric argument shows $t_1$ is well defined.
\end{proof} 

\begin{clm}
The value of $t_0$ is less than $1-\epsilon_0$. 
\end{clm}

\begin{proof}
Suppose $t_0 = 1-\epsilon_0$.  Let $\epsilon> \epsilon_0$ be small enough so that $1 - \epsilon$ is greater than the height of the highest critical submanifold.  As $t_0 = 1-\epsilon_0$ there is a curve $\alpha$ of $F(1 - \epsilon) \cap S$ which is essential in $F(1 - \epsilon)$ but bounds a compression in $H(1 - \epsilon)$.

Recall that the boundary of $S$ has been isotoped into standard position. It follows that the components of $S \cap H'(1 - \epsilon)$ are all disks and bigons. Hence, $\alpha$ bounds compressions for $F(1 - \epsilon)$ on both sides and $d(K,F)=0$.
\end{proof}

\begin{clm}
\label{c:t_0=t_1}
If $t_0=t_1<1-\epsilon_0$ then $d(K,F) = 1$. 
\end{clm}

\begin{proof}
If $t_0=t_1<1-\epsilon_0$ then for all sufficiently small $\epsilon$ there is a curve of $S \cap F(t_0 + \epsilon)$ which bounds a compression in $H'(t)$ and a curve of $S \cap F(t_0 - \epsilon)$ which bounds a compression in $H(t)$. But for $\epsilon$ sufficiently small the curves of $\pi (S \cap F(t_0 + \epsilon))$ can be made disjoint from the curves of $\pi (S \cap F(t_0 - \epsilon))$ because $F$ and $S$ are orientable. This is basically identical to Lemma 4.4 of \cite{gabai:87}. 
\end{proof}

Henceforth we assume that $\epsilon_0< t_0 <t_1<1-\epsilon_0$. 

\begin{clm}
\label{c:atmostone}
If $t_* \in (t_0,t_1)$ is a critical value then for sufficiently small $\epsilon$ the curves of $\pi(F(t_*-\epsilon) \cap S)$ are at a distance of at most one from the curves of $\pi(F(t_*+\epsilon) \cap S)$.
\end{clm}

\begin{proof}
As in the proof of Claim \ref{c:t_0=t_1} for $\epsilon$ sufficiently small the curves of $\pi (S \cap F(t_* + \epsilon))$ can be made disjoint from the curves of $\pi (S \cap F(t_* - \epsilon))$. The result follows unless either of these are collections of inessential curves, and hence are not represented in $\Gamma (F_K)$. However, if this is the case then all curves of $S \cap F(t_* + \epsilon)$ (say) are inessential on $S$. By Lemma \ref{l:EssentialPersists} a sequence of surgeries produces an essential surface disjoint from $F(t_* + \epsilon)$, contradicting Claim \ref{c:noessentialinhandlebody}.
\end{proof}

\begin{clm}
\label{l:OnlySpecialAndMutual}
If a component of $F(t) \cap S$ is inessential on $F(t)$ then it must be inessential on $S$.
\end{clm}

\begin{proof}
This follows directly from Lemma \ref{l:CutImpliesCompression}. 
\end{proof}

Now let $t \in [\epsilon_0, 1-\epsilon_0]$ be a regular value of $h|_S$. Pick a component $\gamma$ of $F(t) \cap S$. The curve $\gamma$ is a {\it mutually essential} curve if it is essential on both $F(t)$ and $S$, {\it mutually inessential} if it is inessential on both and {\it mutual} if it is mutually essential or mutually inessential. Finally, $\gamma$ is {\it special} if it is inessential on $S$, but essential on $F(t)$. Note that there are three kinds of special curves: loops that bound disks on $S$, loops that cobound (with $\partial S$) annuli in $S$, and arcs isotopic (via bigons) into $\partial S$.  

\begin{clm}
\label{c:special}
Suppose $t$ is a regular value of $h|_S$ in $[t_0,t_1]$. Then every curve of $F(t) \cap S$ is mutual.
\end{clm}

\begin{proof}
Pick a regular value $t \in [\epsilon_0,1-\epsilon_0]$. By Claim \ref{l:OnlySpecialAndMutual} we may assume that there is a special curve $\gamma$ in $F(t) \cap S$. By definition, $\gamma$ is essential on $F(t)$ but inessential on $S$. It follows that a component $E$ of $S-\gamma$ is a cut surface. By Lemma \ref{l:CutImpliesCompression} there is a curve of $E \cap F(t)$ which bounds a compression for $F(t)$. This compression either lies in $H(t)$ or $H'(t)$. Since $E \cap F(t) \subset S \cap F(t)$ we conclude $t \notin [t_0,t_1]$. 
\end{proof}

\begin{clm}
\label{c:essentialarcs}
If $\alpha$ is an arc component of $F(t) \cap S$ and $h(\alpha)=t \in (t_0,t_1)$ then $\alpha$ is mutually essential.
\end{clm}

\begin{proof}
By Claim \ref{c:special} the only other possibility is that $\alpha$ is mutually inessential. In this case $\partial \alpha$ is the boundary of some arc $\gamma$ of $\partial S - F(t)$. Also, $\partial \gamma=\partial \alpha$ lies on the same component of $\partial F(t)$. This violates our assumption that $S$ is in standard position. 
\end{proof}

In $h^{-1}([t_0,t_1])$ we see the usual four types of critical submanifolds for $S$: maxima, minima, saddles, and meridional boundary components. Suppose a critical submanifold happens at height $t$ which is a saddle or meridional boundary component. Let $P$ be the closure of the component of $S-F(t \pm \epsilon)$ that contains the critical submanifold. We call $P$ a {\it horizontal neighborhood} (in $S$) of the critical submanifold. Let $\partial _{\pm} P=P \cap F(t \pm \epsilon)$. Then we say the critical submanifold at $t$ is {\it special} if there is some component of $\partial _{\pm} P$ that is special. If the critical submanifold at $t$ is not special then we say it is {\it inessential} if some component of the closure of $S-P$ is a disk and {\it essential} otherwise. Note that if the critical submanifold at $t$ is inessential then it follows from Claim \ref{c:essentialarcs} that there is a mutually inessential loop component of $\partial _{\pm} P$ which bounds a disk in $S$.

\begin{clm}
\label{c:specialsaddle}
Suppose $t_* \in [t_0,t_1]$. If there is a special critical submanifold at $t_*$ then $t_*=t_0$ or $t_1$.
\end{clm}

\begin{proof}
By definition, if a special critical submanifold happens at $t_*$ then there is a special curve $\alpha$ in $S \cap F(t_*-\epsilon)$ or $S \cap F(t_*+\epsilon)$. Assuming the former, Claim \ref{c:special} implies $t_*-\epsilon \notin [t_0,t_1]$. Hence $t_*=t_0$. If, on the other hand, $\alpha \subset F(t_*+\epsilon)$ then we deduce $t_*=t_1$. 
\end{proof}

\begin{lem}
\label{l:LackOfMotion}
Let $t_-$ and $t_+$ be regular values in $[t_0,t_1]$ such that every saddle and every meridional boundary component of $S$ in $h^{-1}(t_-,t_+)$ is inessential. Then $\pi (F(t_-) \cap S)$ and $\pi (F(t_+) \cap S)$ share a vertex in $\Gamma (F_K)$.
\end{lem}

\begin{proof}
Let $\{t_i\}$ be the critical values of $h|_S$ which lie in $[t_-, t_+]$.  Choose $r_i$ slightly greater than the $t_i$
and let $R = \{r_i\} \cup \{t_- + \epsilon\}$.  

For every $r \in R$ surger $S$ in the following way. If $S \cap F(r)$ contains mutually inessential curves then some such curve bounds a cut surface in $F(r)$. Surger $S$ along this cut surface. After a sequence of such surgeries we obtain from $S$ a surface which meets $F(r)$ only in mutually essential curves, for all $r \in R$. 

Let $M'=h^{-1}([t_-,t_+])$. Let $S'$ be the intersection of the surgered surface with $M'$.  Note that $h|_{S'}$, the height function restricted to $S'$, has either two or four new critical values for every surgery performed.  See Figure~\ref{TrivialSurgery}.

\begin{figure}
$$\begin{array}{c}
\epsfig{file=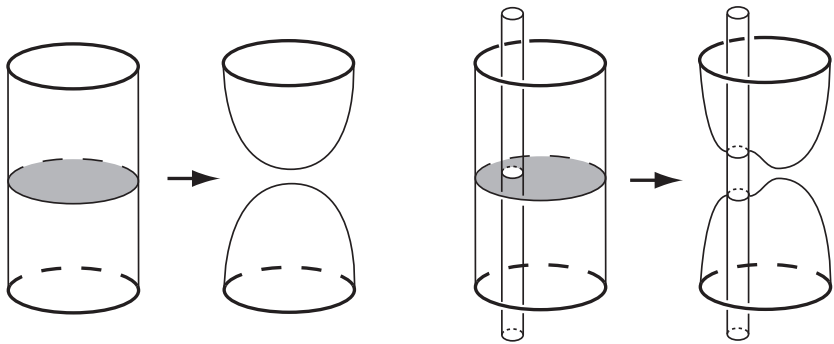, height = 4.5 cm}
\end{array}$$
\caption{Constructing $S'$ from $S$. On the left two new critical values are created. On the right four are created.}
\label{TrivialSurgery}
\end{figure}

We will say a surface $V$ is {\it vertical} in $M'$ if $V=\pi^{-1} (\alpha) \cap M'$, where $\alpha$ is a properly embedded one-manifold in $F_K$. Note that a vertical surface $V$ is either a disk or an annulus.  We need the following claim to prove the lemma:

\begin{claim}
Each component $S''$ of $S'$ is either
\begin{itemize}
    \item a sphere or a meridional annulus, or
    \item properly isotopic into $F(t_-)$ or $F(t_+)$, or
    \item properly isotopic to a vertical surface $V$ with $\pi(V)$ essential in $F_K$.
\end{itemize}
\end{claim}

\begin{proof}
If $h|_{S''}$ has no critical values then $S''$ is isotopic to a vertical annulus or disk.  In this case $S'' \cap \partial M'$ must be essential by the construction of $S'$.  Note that this kind of situation is the desired conclusion of the lemma at hand. If $h|_{S''}$ has only critical values of even index (and no meridional boundary components) then $S''$ is a boundary parallel disk or a sphere.

We now assume that $S''$ contains a critical submanifold which is not a max or min. The component $S''$ either contains a saddle or meridional boundary component of $S$, or it does not. Suppose the latter. It follows that $S''$ is either a meridional annulus or a boundary parallel annulus (with one boundary component on $\partial M_K$). 

Now suppose that $S''$ contains a saddle or meridional boundary component of $S$ at height $t_*$. Let $P$ be the closure of the component of $S''-F(t_* \pm \epsilon)$ that contains this critical submanifold. (Note that $P$ is also a subsurface of $S$ since $\epsilon$ is very small.) Recall that $P$ is the horizontal neighborhood of the critical submanifold. Let $\partial _{\pm} P=P \cap F(t_* \pm \epsilon)$. As every critical submanifold of $S \cap M'$ is inessential at least one loop component of $\partial _{\pm} P$ bounds a disk in $S$ (see the comment which precedes Claim \ref{c:specialsaddle}). 

Now suppose that $S''$ contains a meridional boundary component of $S$ at height $t_*$. Let $P$ be the corresponding horizontal neighborhood. Let $\partial _{\pm} P=C_1 \cup C_2$, where $C_1$ bounds a disk $D$ in $S$. Hence, $D \cup P \subset S$ is a cut annulus and we see that $C_2$ is also inessential in $S$. By Claim \ref{c:special} the $C_i$ are inessential in $F(t_* \pm \epsilon)$. It now follows from Lemma \ref{l:NoTypeMismatch} that $C_1$ bounds a disk in $F(t_* \pm \epsilon)$ while $C_2$ bounds a cut annulus in $F(t_* \mp \epsilon)$. Thus $S''$ is a meridional annulus.

We now assume that $S''$ contains no meridional boundary components of $S$, and hence contains a saddle. Suppose some such saddle has a horizontal neighborhood $P$ such that two components of $\partial _{\pm} P$ are inessential. Then it follows that all three components are inessential. If two bound disks then all three do. Therefore, by Lemma \ref{l:NoTypeMismatch}, $S''$ is a sphere. If one bounds a disk and the other two bound cut annuli then $S''$ is a meridional annulus.

\begin{figure}
\psfrag{x}{\small$x$}
\psfrag{g}{$\gamma _x$}
%\psfrag{c}{\small$F(\theta)$}
%\psfrag{e}{\small$F(\theta-\epsilon)$}
%\psfrag{r}{\small$F(r)$}
$$\begin{array}{c}
\epsfig{file=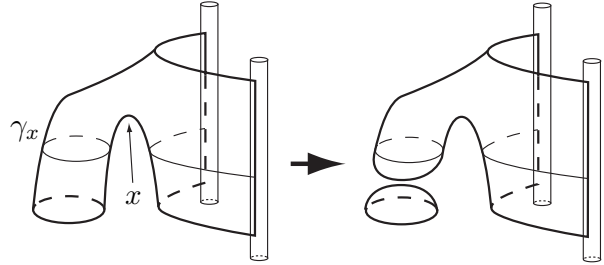, height = 3.5 cm}
\end{array}$$
\caption{Surgery near a saddle whose horizontal neighborhood has exactly one inessential boundary component.}
\label{f:c3surgery}
\end{figure}

Finally, we assume that $S''$ contains no meridional boundary components, and every saddle $x$ has a horizontal neighborhood $P_x$ with exactly one component $\gamma_x$ of $\partial _{\pm} P_x$ inessential, bounding a disk in $S$ (see Figure \ref{f:c3surgery}). By Claim \ref{c:special} and Lemma \ref{l:NoTypeMismatch} it follows that $\gamma_x$ bounds a disk in $S''$. Hence $S''$ is either a union of disks or a union of annuli. In the first case $S''$ is isotopic to a vertical disk. In the latter case $S''$ is either isotopic to a vertical annulus or is a boundary parallel annulus. \end{proof}

Now to complete the proof of the Lemma \ref{l:LackOfMotion}. Suppose that no component of $S'$ meets both boundary components of $M'$.  Thus, by the previous claim, every component of $S'$ meeting $F(t_-)$ is boundary parallel in $M'$. Isotope $F(t_-)$ across these boundary-parallelisms to obtain a surface $F'$ which intersects the surface $S$ only in mutually inessential curves. Some component of $F'-S$ is then a cut surface, which we may use to surger $S$. By Lemma \ref{l:EssentialPersists} we obtain an essential surface which meets $F'$ in fewer curves. Continuing in this fashion we obtain an essential surface disjoint from $F'$, violating Claim \ref{c:noessentialinhandlebody}.

We conclude that there is a component $S'' \subset S'$ which meets both $F(t_-)$ and $F(t_+)$.  By the claim above this $S''$ must be isotopic to a vertical annulus or vertical disk with essential boundary. The lemma is thus proved.
\end{proof}

We now complete the proof of Theorem \ref{t:main}. Note that when $t \in [t_0,t_1]$ is a regular value, $\pi(F(t) \cap S)$ is a properly embedded 1-manifold in $F_K$ (recall that $F_K=F \cap M_K$). The distance between the loops and arcs of $\pi(F(t_0-\epsilon) \cap S)$ and of $\pi(F(t_1+\epsilon) \cap S)$ in $\Gamma (F_K)$ is an upper bound for the distance $d(K,F)$. By Lemma \ref{l:LackOfMotion} and Claim \ref{c:atmostone} this number is bounded by the number of essential critical submanifolds, $e$, plus the number of special critical submanifolds. By Claim \ref{c:specialsaddle} this latter number is at most two. We therefore conclude $d(K,F) \le e+2$.

We now bound the Euler characteristic of $S$. Suppose an essential critical submanifold happens at $t_*$ and let $P$ be its horizontal neighborhood in $S$. Note that in all cases  $\chi(P)=-1$. (When $P$ has vertical boundary compute its Euler characteristic by doubling across the vertical boundary and taking half of the Euler characteristic of the resulting surface. See, for example, the surface on the left in Figure \ref{f:c3surgery}.) By the definition of an essential critical submanifold $\partial P-\partial S$ is essential in $S$. We conclude that $\chi(S)\le -e$.

Putting the above together we conclude:

\begin{eqnarray*}
d(K,F) & \le & e+2\\
& \le & -\chi(S) +2\\
& = & -(2-2g(S)-|\partial S|)+2\\
& = & 2g(S)+|\partial S| 
\end{eqnarray*}

\section{Applications}

We now present a few quick corollaries to Theorem \ref{t:main}. 

\begin{cor}
\label{c:c1}
Suppose $K$ is a knot in $S^3$ whose distance is $d(K,F)$ with respect to a bridge sphere $F$. Then the genus of $K$ is at least $\frac {d(K,F)-1}{2}$.
\end{cor}

\begin{proof}
The genus of $K$ is defined to be the smallest genus of all orientable spanning surfaces for $K$. Such a spanning surface is essential and has exactly one boundary component. Hence, an immediate application of Theorem \ref{t:main} implies $d(K,F) \le 2g(K)+1$. 
\end{proof}

\begin{cor}
\label{c:c2}
If $K$ is a knot whose distance is at least 3 with respect to some Heegaard surface then the complement of $K$ is hyperbolic of finite volume. 
\end{cor}

\begin{proof}
If the distance is greater than two then $M_K$ is irreducible, atoroidal, anannular, and has incompressible boundary. It follows from Thurston's geometrization theorem for Haken manifolds that $M_K$ is hyperbolic of finite volume. 
\end{proof}

For the next corollary we will need the following definition. 

\begin{dfn}
Suppose $M$ is obtained by removing a neighborhood of a knot $K$ in $S^3$ and gluing in a new solid torus to the resulting boundary component. Then we say that $M$ was obtained by {\it Dehn surgery} on $K$. 
\end{dfn}

\begin{cor}
\label{c:c3}
Suppose $K$ is a knot in $S^3$ whose distance is $d(K,F)$ with respect to a bridge sphere $F$. If a manifold $M$, obtained by Dehn surgery on $K$, contains an incompressible torus $T$ then $|\partial (T \cap M_K)|$ is at least $d(K,F)-2$.
\end{cor}
%S^3? T^2?

\begin{proof}
Choose $T$ to minimize $|T \cap K|$ in $M$. Let $T_K=T \cap M_K$. It follows from the minimality assumption that $T_K$ is essential. Theorem \ref{t:main} says that $d(K,F)$ is bounded above by twice the genus of $T_K$ plus $|\partial T_K|$. But $T$ is a torus, so the genus of $T_K$ is one.
\end{proof}

\begin{cor}
\label{c:thinbridge}
Suppose $K$ is a knot in $S^3$ whose distance with respect to some bridge sphere is greater than its bridge number. Then a minimal bridge presentation for $K$ is thin.
\end{cor}

\begin{proof}
Let $F$ be a bridge sphere for which $d(K,F)\ge |K \cap F|$. If thin position for $K$ does not equal bridge position then by \cite{thompson:97} there is a planar, meridional, essential surface $S$ in the complement of $K$ with fewer boundary components than $|K \cap F|$. Hence, by Theorem \ref{t:main} the distance $d(K,F)$ is at most $|\partial S| \le |K \cap F|$.
\end{proof}

\bibliographystyle{alpha}
%\bibliography{master}

\end{document}